\documentclass[12pt]{article} 

\usepackage{times}
\usepackage{amssymb}
\usepackage{amsmath}
\usepackage{graphicx}
\usepackage{epsfig}
\usepackage{subfigure}

\def\shadowbox{\hbox{\rule[-0.0ex]{0.1ex}{1.2ex}%
\hspace{-0.1ex}\rule[-0.0ex]{1.2ex}{0.1ex}%
\hspace{0.0ex}\rule[-0.0ex]{0.1ex}{1.2ex}\hspace{-1.3ex}%
\rule[1.15ex]{1.25ex}{0.1ex}\hspace{-0.0ex}\rule[-0.25ex]{0.3ex}{1.1ex}%
\hspace{-1.2ex}\rule[-0.25ex]{1.1ex}{0.25ex}}}
\def\qed{\ifmmode \hbox{\hfill\shadowbox}
     \else \hphantom{x}\hfill\shadowbox \fi}

\newtheorem{theorem}{Theorem}[section]
\newtheorem{lemma}[theorem]{Lemma}
\newtheorem{definition}[theorem]{Definition}
\newtheorem{proposition}[theorem]{Proposition}
\newtheorem{corollary}[theorem]{Corollary}

\newtheorem{remark}[theorem]{Remark}

\def\remark{{\medskip \noindent \bf Remark:\hspace{0.5em}}}

\def\R {\mathbb{R}}

\def\Nat {\mathbb{N}}
\def\Rtd {\mathbb{R}^{2d}}
\def\I {\mathbb{I}}

\def\K {\mathbb{K}}
\def\G {\mathbb{G}}
\def\Z {\mathbb{Z}}
\def\Zd {\mathbb{Z}^{d}}
\def\Ztd {\mathbb{Z}^{2d}}

\def\Gh {\hat{\mathbb{G}}}
\def\GGh {\mathbb{G}\times\hat{\mathbb{G}}}
\def\GhG {\hat{\mathbb{G}}\times\mathbb{G}}
\def\GGhT {\mathbb{G}\times\hat{\mathbb{G}}\times\mathbb{T}}
\def\D {\mathbb{D}}
\def\B {\mathbf{\mathcal{B}}}
\def\A {\mathbf{\mathcal{A}}}
\def\a {\mathbf{a}}
\def\b {\mathbf{b}}
\def\L2R {L^{2}(\mathbb{R})}

\def\Ls {L_{\sigma}}
\def\s {\hat{\sigma}}
\def\t {\hat{\tau}}
\def\g {\hat{\gamma}}
\def\T {\mathbb{T}}
\def\H {\mathbb{H}}
\def\Nbar {\overline{N}}
\def\h {\mathbf{h}}
\def\hn {\mathbf{h}_{0}}
\def\ho {\mathbf{h}_{1}}

\def\htw {\mathbf{h}_{2}}
\def\i {\mathbf{i}}

\def\e {\mathbf{e}}
\def\ao {\alpha_{1}}
\def\at {\alpha_{2}}
\def\ath {\alpha_{3}}
\def\span {\textnormal{span}}
\def\NKernel {\mathbf{N}^{1}_{v}(\G)}
\def\NovG {\mathcal{N}^{1}_{v}(\G)}
\def\NooG {\mathcal{N}^{1}_{1}(\G)}
\def\Nivg {\mathcal{N}^{\infty}_{v}(\G)}
\def\NovH {\mathcal{N}^{1}_{v}(\H)}
\def\NivH {\mathcal{N}^{\infty}_{v}(\H)}
\def\NooGt {\widetilde{\mathcal{N}^{1}_{1}}(\G)}

\def\NovHt {\widetilde{\mathcal{N}^{1}_{v}}(\H)}
\def\NovGt {\widetilde{\mathcal{N}^{1}_{v}}(\G)}
\def\NivHt {\widetilde{\mathcal{N}^{\infty}_{v}}(\H)}

\def\Pj {P_{j}}
\def\Pi {P_{i}}
\def\Pk {P_{k}}

\def\LovG {L^{1}_{v}(\G)}

\def\LpH {L^{p}(\H)}

\def\WooH {W(L^{1}(\H),l^{1}_{v})}
\def\WioH {W(L^{\infty}(\H),l^{1}_{v})}
\def\LovH {L^{1}_{v}(\H)}

\def\LovGGh {L^{1}_{v}(\GGh)}

\def\hstar {\star}

\begin{document}

\title{\bf Inverse-Closedness of a Banach Algebra of Integral Operators on the Heisenberg Group}

\author{Brendan Farrell and Thomas Strohmer}

\date{\today}
\date{}
\maketitle
\begin{abstract}
Let $\mathbb{H}$ be the general, reduced Heisenberg group. 
Our main result establishes the inverse-closedness of a class of integral operators acting on $L^{p}(\mathbb{H})$, 
given by the off-diagonal decay of the kernel. 
As a consequence of this result, we show that 
if $\alpha_{1}I+S_{f}$, where $S_{f}$ is the operator given by convolution with $f$, $f\in L^{1}_{v}(\mathbb{H})$, 
is invertible in $\B(L^{p}(\mathbb{H}))$, 
then $(\alpha_{1}I+S_{f})^{-1}=\alpha_{2}I+S_{g}$, and $g\in L^{1}_{v}(\mathbb{H})$. 
We prove analogous results for twisted convolution operators and apply the latter results to a class of Weyl pseudodifferential operators.
We briefly discuss relevance to mobile communications.

\vspace{.2in}
\noindent Mathematics Subject Classification: 43A20, 47G30

\vspace{.2in}
\noindent Keywords: Heisenberg group, integral operator, twisted convolution, pseudodifferential operator
\end{abstract}


\begin{section}{Introduction} \label{s:intro}
We consider a class of integral operators defined for the general reduced Heisenberg group $\H$. 
We show that if the kernel of the integral operator $N_{1}$ has $L^{1}_{v}$-integrable off-diagonal decay (here $v$ is a weight) 
and the operator $\alpha_{1} I+N_{1}$ is invertible, then its inverse also has the form $\alpha_{2} I+N_{2}$, 
and the off-diagonal decay is preserved in the kernel of $N_{2}$.
As a consequence of this result, we establish the inverse-closedness of a Banach algebra of convolution operators on 
the general reduced Heisenberg group.
Namely, we consider the operator $S_{f}$ given by convolution with $f$ and
show that if $\alpha_{1}I+S_{f}$, $f\in L^{1}_{v}(\H)$, is invertible in $\B(L^{p}(\H))$, 
then $(\alpha_{1}I+S_{f})^{-1}=\alpha_{2}I+S_{g}$, where $g\in L^{1}_{v}(\H)$. 
While this result relies on other recent results, it has its roots in Wiener's lemma. 

Wiener's lemma states that if a periodic function $f$ has an absolutely summable Fourier series 
\begin{equation}
f(t)=\sum_{n\in\Z}a_{n}e^{2\pi int}\label{series}
\end{equation}
and is nowhere zero, then $1/f$ also has an absolutely convergent Fourier series \cite{Wie32}. 
We call this inverse-closed property the spectral algebra property.
That is, let $\A$ and $\B$, 
$\B\subset\A$, be Banach algebras. 
$\B$ has the \emph{spectral algebra property} if whenever $\b\in\B$ is invertible in $\A$, 
$\b^{-1}\in\B$. 
We equivalently say that $\B$ is \emph{inverse-closed} in $\A$. 
In Wiener's lemma and in many of its descendants, the space $\B$ is given by some form of $l^{1}$ decay, for example 
summability of Fourier coefficients. 
Another perspective on inverse-closedness for a Banach algebra 
views the elements of the algebra as operators and considers to what degree the operator 
maps subspaces to other subspaces or to what degree it leaves subspaces invariant. 
This perspective becomes fruitful when we consider subspaces derived from the structure of the underlying group, 
rather than from a basis. 

To view a periodic function with summable Fourier series $\a=\{a_{n}\}_{n\in\Z}$ as an operator, we consider 
the algebra of bi-infinite Toeplitz matrices with summable antidiagonals, and we let $T_{\a}$ be the Toeplitz matrix given by 
the sequence $\a$ as antidiagonal. 
To frame the Toeplitz case in the subspace perspective, we reformulate the summable antidiagonal property as 
\begin{equation}
\sum_{n\in\Z}\sup_{i-j=n}|\langle T_{\a}e_{i},e_{j}\rangle|=\sum_{n\in\Z}|\langle T_{\a}e_{n},e_{0}\rangle|<\infty,\label{t}
\end{equation}
where $\{e_{i}\}_{i\in\Z}$ is the standard basis for $l^{2}(\Z)$. 
In this case, the subspace perspective simply states that for the subspaces $E_{i}$ and $E_{j}$ 
given by $e_{i}$ and $e_{j}$, $\| T_{a}:E_{i}\rightarrow E_{j}\|\rightarrow 0$ 
as $|i-j|\rightarrow \infty$ with decay given by (\ref{t}).

Bochner and Phillips contributed the first essential step towards a general operator version 
of the spectral algebra property \cite{BP42}. 
They showed that the $a_{n}$ in (\ref{series}) may belong to a--possibly noncommutative--Banach algebra. 
This key result enabled Gohberg, Kaashoek and Woerdeman \cite{GKW89} and Baskakov \cite{Bas90} to establish 
an operator version of the spectral algebra property. 
They considered subspaces $X_{i}$ of the space $X$, indexed by a discrete abelian group $\I$,
satisfying $X_{i}\cap X_{j}=\{0\}$ for $i\neq j$ and 
$X=\overline{\span}\{X_{i}\}_{i\in \I}$, and set $\Pi$ to be the projection onto $X_{i}$. 
For the linear operator $T:\;X\rightarrow X$ they set 
\[a_{n}=\sum_{i-j=n}\Pi T\Pj,\]
and consider the operator-valued Fourier series 
\begin{equation}
 f(t)=\sum_{n\in\I}a_{n}e^{2\pi int}\label{opseries}.
\end{equation}
They then use Bochner and Phillips's work to establish that operators of the form (\ref{opseries}) 
satisfying $\sum_{n\in\Z}\|a_{n}\|=\sum_{n\in\Z}\sup_{i-j=n}\|\Pi T\Pj\|<\infty$ 
form an inverse-closed Banach algebra in $\B(X)$
\cite{GKW89,Bas90,Bas92}. 

In the commutative setting, Gelfand, Raikov and Shilov \cite{GRS64} addressed the important question: 
what rates of decay of an element are preserved in its inverse? 
They answered this question by determining conditions on a weight function $v$ such 
that series finite in the following 
weighted norm form an inverse-closed Banach algebra in $l^{1}_{1}(\Z)$: 
\begin{equation*}
\|\a\|_{l^{1}_{v}(\Z)}=\sum_{n\in\Z}\|a_{n}\|v(n)<\infty.
\end{equation*}
These three conditions on $v$ are given later; the key condition is called the GRS condition, and a function satisfying all three 
is called \emph{admissible}.
Baskakov incorporated the GRS condition and proved the following operator
version of the spectral algebra property\cite{Bas92,Bas97}\footnote{The
version presented here is slightly different from the theorems
in~\cite{Bas92,Bas97}, but it can be easily extracted from the proof
of Theorem~2 in~\cite{Bas92}.}:
let $v$ be an admissible weight; if the linear operator $T$ satisfies 
\begin{equation*}
\sum_{n\in\Z}\sup_{i-j=n}\|\Pi T\Pj\|v(n)<\infty
\end{equation*}
and is invertible, then 
\begin{equation*}\sum_{n\in\Z}\sup_{i-j=n}\|\Pi T^{-1}\Pj\|v(n)<\infty.
\end{equation*} 
Kurbatov considered a class of operators satisfying 
\begin{equation}(Tf)(t)\leq\int\beta(t-s)|f(s)|ds\label{maj}\end{equation}
for some $\beta\in L^{1}$.
He showed, using results very similar to Baskakov's (but derived
independently), that if $\alpha_{1}I+T_{1}$ is invertible and $T_{1}$ satisfies (\ref{maj}) 
for $\beta_{1}\in L^{1}$, then $(\alpha_{1}I+T_{1})^{-1}=\at I+T_{2}$ and $T_{2}$ 
satisfies $(\ref{maj})$ for $\beta_{2}\in L^{1}$ \cite{Kur01}. 
This theorem, as stated for integral operators in \cite{Kur99}, is the point of departure for the research presented in this paper.

We have two motivations for the work presented here: on the one hand a 
question of abstract harmonic analysis and on the other hand research on 
the propagation channel of a mobile communication system. 
The abstract harmonic analysis question is: for what nonabelian groups does the spectral algebra property hold? 
One recent result in this direction is by Gr\"ochenig and Leinert \cite{GL04,GL06}. 
They established the spectral algebra property for $(l^{1}_{v}(\Ztd\times\Ztd),\natural_{\theta})$,
where 
$\natural_{\theta}$ is the following form of twisted convolution:
\begin{equation*}
(\a\natural_{\theta}\b)_{(m,n)}=\sum_{k,l\in\Zd}a_{kl}b_{m-k,n-l}e^{2\pi i\theta(m-k)\cdot l},
\end{equation*}
and they use this result to prove the spectral algebra property for convergent sums of time-frequency shifts 
$\sum_{\lambda\in\Lambda}c_{\lambda}T_{x_{\lambda}}M_{\omega_{\lambda}}$, $\sum_{\lambda\in\Lambda}|c_{\lambda}|v(\lambda)<\infty$. 
($T_{x}$ and $M_{\omega}$ are defined below.) 
Balan recently generalized this result by relaxing the lattice requirement to solely a discrete \emph{subset} of $\Rtd$. 
He showed that if $\sum_{\lambda\in\Lambda}c_{\lambda}T_{x_{\lambda}}M_{\omega_{\lambda}}$,  
$\sum_{\lambda\in\Lambda}|c_{\lambda}|v(\lambda)<\infty$, 
is invertible, then 
$(\sum_{\lambda\in\Lambda}c_{\lambda}T_{x_{\lambda}}M_{\omega_{\lambda}})^{-1}=\sum_{\sigma\in\Sigma}c_{\sigma}T_{x_{\sigma}}M_{\omega_{\sigma}}$ and $\sum_{\sigma\in\Sigma}|c_{\sigma}|v(\sigma)<\infty$,
where $\Lambda$, $\Sigma\subset\Rtd$, and $|\Lambda|$, $|\Sigma|<\infty$, but $\Lambda\neq\Sigma$ (in general) \cite{Bal07}.

Many of the other recent results of this nature are for integral operators, where the spectral algebra property 
is manifested in the kernel. 
This is true for the seminal paper by Sj\"ostrand in pseudodifferential operator theory \cite{Sjo95}. 
In that paper he proves the matrix version of Baskakov's result, and uses this to prove the spectral algebra 
property for pseudodifferential operators with symbols in $M^{\infty,1}_{v}(\Rtd)$. 
(Let $g\in\mathcal{S}(\Rtd)$ be a compactly supported, $C^{\infty}$ function 
satisfying $\sum_{k\in\Z}g(t-k)=1$ for all $t\in\Rtd$. Then the symbol $\sigma\in\mathcal{S}'(\Rtd)$ belongs 
to $M^{\infty,1}_{v}(\Rtd)$ if $\int_{\Rtd}\sup_{k\in\Ztd}|(\sigma\cdot g(.-k))^{\wedge}(\zeta)|d\zeta<\infty$.)
Sj\"ostrand proved the matrix Baskakov result and his pseudodifferential operator result using techniques from 
``hard analysis''. 
Gr\"ochenig later also achieved this same result using techniques solely from harmonic analysis \cite{Gro04a}. 
In fact, he proved more. He showed that a pseudodifferential operator $L_{\sigma}$ has Weyl symbol 
$\sigma\in M^{\infty,1}_{v}(\Rtd)$ if and only if the matrix given by 
$M_{m,n,m',n'}=\langle L_{\sigma}\psi_{m,n},\psi_{m',n'}\rangle$, 
for a proper Gabor frame $\{\psi_{k,l}\}_{k,l\in\Z}$ \cite{Gro01}, is in the Baskakov matrix algebra 
given by $\sum_{(k,l)\in\Z^{2}}\sup_{(m-m',n-n')=(k,l)}|A_{m,n,m',n'}|v(k,l)<\infty$ \cite{Gro04a}. 
Gr\"ochenig and one of us recently extended this result to pseudodifferential operators with symbols 
defined on $\GGh$, where $\G$ is any locally compact abelian group and $\Gh$ is its dual group \cite{GS06}. 
Locally compact abelian groups, their dual groups and twisted convolution are standard features throughout the work 
just discussed; therefore, it is very natural to look directly at the general reduced Heisenberg group for a fundamental theorem.

In mobile communications a transmitted signal travels through a channel that is modeled by a pseudodifferential operator. 
When a single source transmits a signal it is reflected by objects in its environment, which results in 
different paths from transmitter to receiver, each with its own travel time. 
In the case of mobile communications, a moving transmitter and/or receiver
gives rise to the Doppler effect~\cite{Rap96}, which results in a frequency shift.
Thus, denoting time shift by $T_{x}f(t)=f(t-x)$ and modulation or frequency shift by $M_{\omega}f(t)=e^{2\pi i\omega t}f(t)$, 
the received signal can be represented as the following collection of weighted, delayed and modulated 
copies of the transmitted signal:
\begin{equation*}
f_{rec}(t)=\int_{\R}\int_{\R^{+}}\s(\omega,x)T_{x}M_{\omega}f_{trans}(t)dxd\omega\label{channel}.
\end{equation*}
We, therefore, consider Weyl pseudodifferential operators:
\begin{equation*}
\Ls f(t)=\int_{\G}\int_{\Gh}\hat{ \sigma}(\omega,x)e^{-\pi i \omega\cdot x}T_{-x}M_{\omega}f(t)d\omega dx.
\end{equation*} 
In particular, we posed the question, if $\alpha_{1}I+\Ls$ is invertible and 
$\s\in L^{1}_{v}(\GhG)$, does $(\alpha_{1}I+\Ls)^{-1}=\alpha_{2}I+L_{\tau}$, where $\t\in L^{1}_{v}(\GhG)$? 

Throughout this paper, $\widetilde{\A}$ will denote the Banach algebra $\A$ with adjoined identity. 
In Section 2 we prove the spectral algebra property for 
$\NovHt$, where $\NovH$ is the space of integral operators with kernels having $\LovH$-integrable off-diagonal decay.
The basis for our proof is establishing an operator class for which we can apply Baskakov's theorem, 
and using a dense, two-sided, proper ideal within that class.  
We use this result to prove the spectral algebra property for convolution operators on $L^{1}_{v}(\H)$.
In Section 3 we prove the spectral algebra property for $(\widetilde{L^{1}_{v}}(\GGh),\natural)$, where $\natural$ is twisted convolution. 
This result is in the same spirit as work of Gr\"ochenig and Leinert on 
$(L^{1}_{v}(\Zd\times\Zd),\natural_{\theta})$ \cite{GL04,GL06}, 
but more general, in that it holds for arbitrary locally compact abelian groups. 
We apply these results to the class of pseudodifferential operator with symbols $\sigma$
satisfying $\s\in L^{1}_{v}(\GGh)$. Lastly, we discuss the consequences of these theorems for mobile communication channels.

\end{section}

\begin{section}{The Convolution Algebra on the General, Reduced Heisenberg Group}

Our construction of the general, reduced Heisenberg group begins with the locally compact abelian group $\G$ 
and its dual group $\Gh$, which is also locally compact and abelian \cite{Fol94}. 
We assume that $\G$ is second countable and metrizable.
Throughout this paper, arbitrary groups will be denoted $G$, and locally compact abelian groups 
will be denoted $\G$. 
While $\H$ is not abelian, it will still be written in the same font as $\G$.
Elements of $\G$ will be written in Latin letters and elements of $\Gh$  in Greek letters. 
By Pontrjagin's duality theorem, $\hat{\Gh}\cong\G$ \cite{RS00}.
For convenience, we will set 
$e^{2\pi ix\cdot\omega}=\langle x,\omega\rangle$, where $\langle x,\omega\rangle$ denotes the action of 
$\omega\in\Gh$ on $x\in\G$. 

Here we approach the Heisenberg group from the perspective of pseudodifferential operators and time-frequency analysis, 
and thus motivate it from the operators translation and modulation. 
Translation is right addition by the inverse of an element in $\G$: $T_{x}f(y)=f(y-x)$; 
modulation is multiplication by the evaluation of a character in $\Gh$: 
$M_{\omega}f(y)=\langle \omega,y\rangle f(y)$.
However the set of operators $T_{x}M_{\omega}$ is not closed, as 
$(T_{x}M_{\omega})(T_{x'}M_{\omega'})=e^{2\pi i x'\cdot\omega}T_{x+x'}M_{\omega+\omega'}$,
and therefore this set of operators is not parameterizable by $\G\times\Gh$,
but by $\G\times\Gh\times\T$. 
The extension of $\GGh$ to $\GGhT$ is called the 
general, reduced Heisenberg group $\H=\GGh\times\T$ \cite{Fol89,Gro01}. 
Elements of $\H$ will be written in bold, and elements of $\G$, $\Gh$ and $\T$ will be written in the normal font.
The group operation for $\H$ is written as multiplication, while the operations for $\G$ and $\Gh$ are written additively 
and for $\T$ is written multiplicatively: 
\begin{equation*}
 \h\h'=(x,\omega,e^{2\pi i\tau})(x',\omega',e^{2\pi i\tau'})=(x+x',\omega+\omega',e^{2\pi i(\tau+\tau')}e^{\pi i(x'\cdot\omega-x\cdot\omega')}).
\end{equation*}
The identity on $\H$ is $\e=(0,0,1)$. 
The measure on $\H$ is $d\h=dxd\omega d\tau$, where $dx,d\omega$ and $d\tau$ correspond to the 
invariant measures on $\G,\Gh$ and $\T$ respectively, normalized so that the measures of $U_{\G},U_{\Gh}$ and 
$U_{\T}$, to be defined shortly, are each 1. 
Since $\G,\Gh$ and $\T$ are commutative, the invariant measure is both left and right invariant.
The space $L^{1}_{v}(\H)$ consists of those functions satisfying 
\begin{equation*}
\|f\|_{L^{1}_{v}(\H)}=\int_{\H}|f(\h)|v(\h)d\h,
\end{equation*}
where $v$ is an admissible weight, as defined below.
We use $\hstar$ to denote the convolution of two functions defined on $\H$ as follows
\begin{equation*}
 (F_{1}\hstar F_{2})(\hn)=\int_{\H}F_{1}(\h)F_{2}(\h^{-1}\hn)d\h.
\end{equation*}
We now address three preliminaries: partitions, weight functions and the amalgam spaces.
\begin{definition}
Let $G$ be a group. 
$(\mathcal{I},U)$ is a \emph{partition} of $G$ if $\mathcal{I}$ is a discrete set, 
$U$ is subset of a locally compact group, 
$(iU)\bigcap (i'U)=\emptyset$ for $i\neq i'$, and
$\bigcup_{i\in\mathcal{I}}(iU)$ covers $G$. 
For simplicity we assume that $U$ contains the identity.
\end{definition}

\begin{lemma}\label{PartH}
The general reduced Heisenberg group $\H$ possesses a partition.
\end{lemma}
\begin{proof}
We call on the structure theorem, which states that
for any locally, compact, abelian group $\G$, 
$\G\cong \R^{d}\times \G_{0}$, where the locally compact abelian group $\G_{0}$ contains a 
compact, open subgroup $\K$ \cite{HR63}. 
$\D=\G_{0}/\K$ is a discrete group and $\G_{0}=\bigcup_{d\in \D}(d\K)$, 
$\G=\bigcup_{(i,d)\in \Z^{d}\times\D}((i,d)([0,1)\times \K))$, 
and different blocks are disjoint.
If the group $\G_{0}$ contains the compact open subgroup $\K$, then $\hat{\G}_{0}$ contains the 
compact open subgroup $\K^{\perp}$ \cite{RS00}; thus a partition  of $\Gh$ exists that is
analogous to the partition used for $\G$.
While the structure theorem does not apply in general to nonabelian groups, 
we can still partition $\H$ into blocks by the following construction. 
Set $\D=(\G_{0}\times\Gh_{0})/(\K\times\K^{\perp})$, 
$\mathcal{I}=\Z^{2d}\times\D\times\{0\}$, and $U=[0,1)^{2d}\times\K\times\K^{\perp}\times\T$.  
Then $\bigcup_{\i\in\mathcal{I}}(\i U)$ covers $\H$ and $(\mathcal{I},U)$ is a partition for $\H$.
Note that $\mathcal{I}$ is not closed, and hence is not a group.
\end{proof}
\begin{definition} Let $G$ be a group. A weight function $v$ defined on $G$ is 
\emph{admissible} if it satisfies the following three conditions:
 \begin{enumerate}
\item $v$ is continuous, symmetric, i.e. $v(x)=v(x^{-1})$, and normalized so that $v(0)=1.$
\item $v$ is submultiplicative, i.e. $v(xy)\leq v(x)v(y)$ for all $x,y\in G$.
\item $v$ satisfies the Gelfand-Raikov-Shilov (GRS) \cite{GRS64} condition:
\[\lim_{n\rightarrow\infty}v(nx)^{1/n}=1\;\;\;\textnormal{for all }x\in G.\]
\end{enumerate}
\end{definition}
\begin{note}
Throughout this paper the weight $v$ will be assumed to be admissible.
\end{note}
\begin{definition}\label{amalgam}
Let $G$ be a locally compact group and $v$ an admissible weight function.
The \emph{amalgam space} $W(L^{p}(G),l^{q}_{v})$ is the space of functions finite in the 
local $L^{p}$ norm and the global $l^{q}_{v}$ norm as follows:
\begin{equation*}
 \|f\|_{W(L^{p}(G),l^{q}_{v})}=\left(\sum_{i\in \mathcal{I}}\|f\|^{q}_{L^{p}(iU)} v(i)^{q}\right)^{1/q},
\end{equation*}
for a partition $(\mathcal{I},U)$ of $G$.
\end{definition}
The definition of $W(L^{p}_{v}(G),l^{q}_{v})$ is independent of the partition, as different partitions result 
in equivalent norms; see \cite{Ste79,Hei03a}. 
Note that $L^{1}_{v}(G)=W(L^{1}(G),l^{1}_{v})$.
Kurbatov uses the amalgam spaces to prove the inverse-closedness of a class of integral operators 
given by the off-diagonal decay of the kernel. 
He considers integral operators $N$ of the form 
\begin{equation*}
(Nf)(t)=\int_{s\in \G}n(t,s)f(s)ds,
\end{equation*}
where $\G$ is a locally compact abelian group. 
We introduce some notation and define three spaces: 
$\NKernel$ is the space of kernels $n$ for which there exists $\beta\in \LovG$ s.t. $|n(t,s)|\leq\beta(ts^{-1})$ for all $t,s\in \G$. 
$\NovG$ is the space of integral operators with kernel $n\in \NKernel$. 
An operator satisfying this property is said to be \emph{majorized} by $\beta$.
For $\H$ we define $Q_{i,d}=(i,d,0)([0,1)^{2d}\times\K\times\K^{\perp}\times\T)$, 
where $(i,d)\in\Z^{2d}\times\D$, as in Lemma~\ref{PartH}.
To make notation easier, $\mu,\nu$ and $\gamma$ will be elements of $\Z^{2d}\times\D$.
Define $P_{\mu}$ to be the projection of $W(L^{p}(\H),l^{q}_{v})$ onto 
$\{f|f\in W(L^{p}(\H),l^{q}_{v}),\textnormal{supp}(f)\subset Q_{\mu}\}$.
For $N\in\B(W(L^{p}(\H),l^{q}_{v}))$, 
we set $N_{\mu,\nu}=P_{\mu}NP_{\nu}$, and $N_{\gamma}=\sum_{\mu\nu^{-1}=\gamma}N_{\mu,\nu}$.
$\Nivg$ is the class of operators (not necessarily integral operators) satisfying 
\begin{equation*}
 \sum_{\gamma\in \Z^{2d}\times\D}\sup_{\mu\nu^{-1}=\gamma}\|N_{\mu,\nu}:L^{1}(Q_{\nu})\rightarrow L^{\infty}(Q_{\mu})\|v(\gamma)<\infty.
\end{equation*}

The following is a slightly restricted version of Kurbatov's theorem as it applies to the work in this paper:
\begin{theorem}\label{Kurbatov} \textnormal{Kurbatov, Theorem 5.4.7 
\cite{Kur99}.}
Let $\G$ be a non-discrete, locally compact abelian group. 
$\NooGt$ is an inverse-closed subalgebra of $\B(W(L^{p}(\G),l^{q}_{v}))$ 
for $1\leq p,q\leq\infty$.
\end{theorem}
$\remark$
If $\G$ is discrete it is not necessary to adjoin the identity operator to $\NooG$.
The corresponding version of Theorem \ref{Kurbatov} for a discrete 
group $\G$ is a special case of Baskakov's more general and very significant Theorem 1 in \cite{Bas97}. 

We show that Kurbatov's result, Theorem~\ref{Kurbatov} above, holds with admissible weight functions for the nonabelian group $\H$. 
\begin{theorem}\label{KurbatovH}
Let $\mathcal{N}^{1}_{v}(\H)$ denote those bounded integral operators $N$ on $W(L^{p}_{v}(\H),l^{q})$, $1\leq p,q\leq\infty$, 
of the form 
\begin{equation*}
 (Nf)(\hn)=\int_{\H}n(\hn,\h)f(\h)d\h,
\end{equation*}
for which there exists $\beta\in L^{1}_{v}(\H)$ satisfying
\begin{equation*} 
|n(\hn,\h)|\leq \beta(\hn^{-1}\h)
\end{equation*}
for all $\hn,\h\in \H$.
Then $\NovHt$ is an inverse-closed Banach algebra in $\B(W(L^{p}_{v}(\H),l^{q}))$.
\end{theorem}

Before we prove Theorem~\ref{KurbatovH} we need some
preparation. Kurbatov's proof of Theorem \ref{Kurbatov} can be adapted 
to the nonabelian group $\H$,
Theorem~\ref{KurbatovH}, 
once two essential pieces are established. 
First, an appropriate partition must be developed for $\LovH$ that allows us to 
apply Baskakov's result. 
Second, one must establish that $\NivH$ is a two-sided ideal in $\NovH$. 
The proofs below of the intermediate results,
Propositions \ref{5.3.4} and \ref{5.4.5} and Theorem \ref{Kurbatov} follow very closely Kurbatov's 
proofs of the analogous results for the abelian case, cf.\ Sections~5.3 
and 5.4 of \cite{Kur99}.

\begin{lemma}\label{NoIdentity} 
The identity operator $I$ is not an element of $\NivH$.
\end{lemma}
\begin{proof}
Let $(\mathcal{I},U)$ be a partition for $\H$. 
Consider the indicator function $\chi_{E}$, where $E\subset iU$ for some $i\in\mathcal{I}$, 
and the measure of $E$ is $\mu(E)=\epsilon>0$. 
Then $\|I:L^{1}(iU)\rightarrow L^{\infty}(iU)\|\geq 1/\epsilon.$
\end{proof}

\begin{lemma} \textnormal{Kurbatov, Lemma 5.3.3 \cite{Kur99}.}\label{Kur5_3_3}
Let $Q$ and $\overline{Q}$ be locally compact topological spaces with measures $\lambda$ and $\overline{\lambda}$, respectively.
For any $N\in \B(L^{1}(Q),L^{\infty}(\overline{Q}))$ there exists a function $n\in L^{\infty}(Q\times\overline{Q})$ 
such that for all $f\in L^{1}(Q)$ one has 
\begin{equation*}
(Nf)(t)=\int n(t,s)x(s)d\lambda(s).
\end{equation*}
\end{lemma}

\begin{proposition} \label{5.3.4}
The operator $N$ is an integral operator majorized by $\WioH$ if and only if 
$\sum_{\gamma}\sup_{\mu\nu^{-1}=\gamma}\| N_{\mu,\nu}:L^{1}(Q_{\nu})\rightarrow L^{\infty}(Q_{\nu})\|v(\gamma)<\infty$. 
That is, $\NivH$ is the class of integral operators with kernels majorized by functions in $\WioH$. 
\end{proposition}
\begin{proof}
If $N$ is majorized by $\beta\in \WioH$, 
then $\| N_{\mu,\nu}:L^{1}(Q_{\nu})\rightarrow L^{\infty}(Q_{\mu})\|\leq\sup_{Q_{\mu\nu^{-1}}}|\beta|$, which proves the first claim.
We prove the second claim.
Let $f_{\mu}=P_{\mu}f$, where $P_{\mu}$ is the projection onto $Q_{\mu}$, as defined following Definition~\ref{amalgam}. 
$N\in \NivH$ implies $(Nf)_{\mu}=\sum_{\nu\in\Ztd\times\D}N_{\mu\nu}f_{\nu}$. 
In the following, $r,s$ and $t$ will be elements of $\R^{2d}\times\K\times\K^{\perp}\times\T$.
By Lemma \ref{Kur5_3_3}, there exists $n_{\mu\nu}\in L^{\infty}(Q_{\mu}\times Q_{\nu})$ such that 
\begin{equation*}
 (N_{\mu\nu}f_{\nu})(t)=\int_{Q_{\nu}}n_{\mu\nu}(t,s)f_{\nu}(s)d\lambda(s)
\end{equation*}
for $t\in Q_{\mu}$. 
Setting 
\begin{equation*} \alpha_{\gamma}=\sup_{\gamma=\mu\nu^{-1}}\|N_{\mu\nu}:L^{1}(Q_{\nu})\rightarrow L^{\infty}(Q_{\mu})\|,
\end{equation*}
$\gamma\in\Ztd\times\D$, we have $|n_{\mu\nu}(t,s)|\leq\alpha_{\mu\nu^{-1}}$ for all $t\in Q_{\mu}$ and $s\in Q_{\nu}$.
Defining $n(t,s)$ to equal $n_{\mu\nu}(t,s)$ for $t\in Q_{\mu}$, $s\in Q_{\nu}$, we have 
\begin{eqnarray*}
 (Nf)(t)&=&\sum_{\mu\in\Ztd\times\D}\sum_{\nu\in\Ztd\times\D}\int_{Q_{\nu}}n_{\mu\nu}(t,s)f(s)d\lambda(s)\\
&=& \int n(t,s)f(s)d\lambda(s).
\end{eqnarray*}
We now must show that $n$ is majorized by a function $\beta\in \WioH$.
\[Q_{\mu}Q_{\nu}^{-1}=(i_{\mu}i_{\nu}^{-1},d_{\mu}d_{\nu}^{-1},0)((-1,1)^{2d}\times\K\times\K^{\perp}\times\T),\]
where $\mu=(i_{\mu},d_{\mu})$ and $\nu=(i_{\nu},d_{\nu})$. 

For $r\in \R^{2d}\times\K\times\T$, we define 
\begin{equation*}\beta(r)=\sup\{\alpha_{\mu}:(r)\in (i_{\mu},d_{\mu},0)((-1,1)^{2d}\times\K\times\K^{\perp}\times\T)\}.
\end{equation*}
We then have $|n(t,s)|\leq\alpha_{\mu\nu^{-1}}\leq\beta(ts^{-1})$ for $t\in Q_{\mu}$ 
and $s\in Q_{\nu}$. 
We define $\Delta=\Delta(r)$ to be the set of all grid points $\mu\in \Z^{2d}\times\D$ 
such that $r\in (i_{\mu},d_{\mu},0)((-1,1)^{2d}\times\K\times\K^{\perp}\times\T)$.
Now $\beta(r)\leq \max\{\alpha_{\mu}:\mu\in\Delta\}$, 
which implies that $\beta(r)\leq\sum_{\mu\in\Delta}\alpha_{\mu}$. 
By the definition of $\mathcal{N}^{\infty}_{v}(\H)$, $\sum_{\mu}\alpha_{\mu}v(\mu)<\infty$. 
$\Delta$ has at most $2^{2d}$ elements for any $r$. 
Since $\beta$ is constant on each block $Q_{\mu}$, $\|\beta\|_{\WioH}\leq 2^{2d}\sum_{\mu}\alpha_{\mu}v(\mu)<\infty$.
\end{proof}
\begin{proposition}\label{5.4.5} 
$\NivH$ is dense in $\NovH$.
\end{proposition}
\begin{proof}
By Proposition \ref{5.3.4}
we may use that $\WioH$ is dense in $\WooH$, and choose $\overline{\beta}\in\WioH$ 
such that $\|\beta-\overline{\beta}\|_{\WooH}<\epsilon$. 
We may assume that $0\leq\overline{\beta}(\h)\leq\beta(\h)$ for all $\h$. 
Set 
\begin{equation*} \overline{n}(\hn,\ho)=\left\{ \begin{array}{lcl}\frac{\overline{\beta}(\hn\ho^{-1})}{\beta(\hn\ho^{-1})}&:&\beta(\hn\ho^{-1})\neq 0\\0&:&\beta(\hn\ho^{-1})=0\end{array}\right. 
\end{equation*}
We then have $\overline{n}(\hn,\ho)\leq\overline{\beta}(\hn\ho^{-1})$ 
and $0\leq n(\hn,\ho)-\overline{n}(\hn,\ho)\leq \beta(\hn\ho^{-1})-\overline{\beta}(\hn\ho^{-1})$.
\end{proof}
The following is one of the three cases covered by Baskakov's Theorem 1 in \cite{Bas92}.
\begin{theorem}\label{BaskakovFull}
Let $v$ be an admissible weight function, $\I$ a discrete abelian group, and $\{X_{i}\}_{i\in\I}$ 
subspaces of $X$ satisfying $X_{i}\cap X_{j}=\{0\}$ for $i\neq j$ and $X=\overline{\span}\{X_{i}\}_{i\in\I}$. 
Let $P_{i}$ be the projection onto $X_{i}$. 
If $T$ is invertible in $\B(X)$ and 
\begin{equation*}
\sum_{i\in\I}\sup_{j-k=i}\|\Pk T\Pj\|v(n)<\infty
\end{equation*}
then 
\begin{equation*}\sum_{i\in\Z}\sup_{j-k=i}\|\Pj T^{-1}\Pk\|v(n)<\infty.
\end{equation*}
\end{theorem} 
\begin{theorem}\label{Baskakov}
If $N\in \NivHt$ is invertible in $\B(W(L^{p}_{v}(\H),l^{q}))$, then $N^{-1}\in \NivHt$.
\end{theorem}
\begin{proof}
We first define the space of operators 
\[M_{v}=\{T|\sum_{\gamma\in\Z^{2d}\times\D}\sup_{\mu\nu^{-1}=\gamma}\|T:L^{p}(Q_{\nu})\rightarrow L^{p}(Q_{\mu})\|v(\gamma)<\infty\;\;\forall p\in[1,\infty]\},\]
where for each $p$, $T$ is understood to be identically defined on the common part of different spaces. 
(Note that $L^{\infty}(Q_{\mu})$ is dense in each $L^{p}(Q_{\mu})$, $p\in[0,\infty]$.)
By Theorem \ref{BaskakovFull}, if $T\in M_{v}$ and $T$ is invertible in $\B(W(L^{p}(\H),l^{q}_{v}))$, then $T^{-1}\in M_{v}$.
$I$ is clearly an element of $M_{v}$, though it is not an element of $\NivH$ by Lemma \ref{NoIdentity}. 
For $N\in\NivH$ and $T\in M_{v}$, 
\begin{eqnarray*}
\lefteqn{ \sum_{\gamma\in\Z^{2d}\times\D}\sup_{\mu\nu^{-1}=\gamma}\|NT:L^{1}(Q_{\nu})\rightarrow L^{\infty}(Q_{\mu})\|v(\gamma)}\\
&\leq& \sup_{\mu,\nu}\|T:L^{\infty}(Q_{\nu})\rightarrow L^{\infty}(Q_{\mu})\| \sum_{\gamma\in\Z^{2d}\times\D}\sup_{\mu\nu^{-1}=\gamma}\|N:L^{1}(Q_{\nu})\rightarrow L^{\infty}(Q_{\mu})\|v(\gamma),
\end{eqnarray*}
and 
\begin{eqnarray*}
\lefteqn{ \sum_{\gamma\in\Z^{2d}\times\D}\sup_{\mu\nu^{-1}=\gamma}\|TN:L^{1}(Q_{\nu})\rightarrow L^{\infty}(Q_{\mu})\|v(\gamma) }\\
&\leq& \sup_{\mu,\nu}\|T:L^{1}(Q_{\nu})\rightarrow L^{1}(Q_{\mu})\| \sum_{\gamma\in\Z^{2d}\times\D}\sup_{\mu\nu^{-1}=\gamma}\|N:L^{1}(Q_{\nu})\rightarrow L^{\infty}(Q_{\mu})\|v(\gamma);
\end{eqnarray*}
therefore, $\NivH$ is a proper two-sided ideal in $M_{v}$. 
If $\ao I+N_{1}\in\NivHt$ and $\ao I+N_{1}$ is invertible in $\B(L^{p}(\H))$, then Theorem \ref{BaskakovFull} implies
$(\ao I+N_{1})^{-1}=T$ for some $T\in M_{v}$. 
Then $(\ao I+N_{1})T=I$ implies $T=\frac{1}{\ao} I-\frac{1}{\ao} N_{1}T$. 
By the ideal property of $\NivH$, $T=\at I+N_{2}$ for some $N_{2}\in \NivH$.

\end{proof}

\begin{lemma}\label{young}
$\NivH$ is a two-sided ideal in $\NovH$.
\end{lemma}
\begin{proof}
Assume $f\in L^{1}_{v}(\H)$ and $g\in L^{\infty}_{v}(\H)$. Then 
\begin{eqnarray*}
\|f\hstar g\|_{L^{\infty}_{v}(\H)}&=&\sup_{\hn}\left|\int_{\H}f(\h)g(\h^{-1}\hn)d\h\right| v(\hn)\\
&\leq&\int_{\H}|f(\h)|\sup_{\hn}|g(\h^{-1}\hn)v(\hn)|d\h\\
&\leq&\int_{\H}|f(\h)|\sup_{\hn}|g(\hn)v(\hn)|v(\h)d\h\\
&=&\|f\|_{ L^{1}_{v}(\H)}\|g\|_{ L^{\infty}_{v}(\H)}.
\end{eqnarray*}
One similarly shows that $L^{\infty}_{v}(\H)\hstar L^{1}_{v}(\H)\subset L^{\infty}_{v}(\H)$.
By Theorem 11.8.3 in \cite{Hei03a} and the discussion immediately following it concerning $\H$,
$\WioH$ is a two-sided ideal in $\WooH$ with respect to convolution. 
The lemma then follows from the  composition rule for majorized integral operators 
given at the start of the proof of Theorem \ref{KurbatovH}. 

\end{proof}
\begin{proof}\textbf{of Theorem \ref{KurbatovH}}
Let $N_{1},\;N_{2}\in\NovH$ be majorized, respectively, by $\beta_{1}$ and $\beta_{2}$. 
Using Fubini's theorem, we have
\begin{eqnarray*}
(N_{1}N_{2})f(\hn)&=&\int n_{1}(\hn,\ho)\int n_{2}(\ho,\htw)f(\htw)d\htw d\ho\\
&=& \iint n_{1}(\hn,\ho)n_{2}(\ho,\htw)f(\htw)d\ho d\htw\\
&=&\int n(\hn,\htw)f(\htw)d\htw.
\end{eqnarray*}
Therefore, $N_{1}N_{2}$ defines an integral operator of the same form.
\begin{eqnarray*}
(N_{1}N_{2})f(\hn)&=&\int n_{1}(\hn,\ho)\int n_{2}(\ho,\htw)f(\htw)d\htw d\ho\\
&\leq & \int  \beta_{1}(\hn^{-1}\ho)\int \beta_{2}(\ho^{-1}\htw)|f(\htw)|d\htw d\ho\\
&=&\iint \beta_{1}(\ho)\beta_{2}(\ho^{-1}\hn^{-1}\htw)d\ho|f(\htw)|d\htw\\
&=&\int\beta(\hn^{-1}\htw)|f(\htw)|d\htw
\end{eqnarray*}
for $\beta=\beta_{1}\hstar\beta_{2}\in L^{1}_{v}(\H)$.
This establishes the Banach algebra property. 

Assume that the operator $\alpha I+N$, $N\in \NovH$, is invertible. 
We first show that $I\notin \NovH$. 
If $I\in\NovH$, then since $\NivH$ is dense in $\NovH$ (Proposition \ref{5.4.5}), 
$\NivH$ would contain a sequence approaching $I$. 
Lemma \ref{NoIdentity} shows that such an operator would be unbounded in the $\NivH$-norm. 
Therefore $\alpha\neq 0$.

By Proposition~\ref{5.4.5} we may choose $\Nbar\in\NivH$ such that 
$\|n-\overline{n}\|_{\WooH}<\alpha/2$. 
Then by Proposition~1.4.2 in~\cite{Kur99}, $\alpha I +(N-\Nbar)$ is invertible in $\NovHt$. 
As in the proof of Theorem 5.4.7 in~\cite{Kur99} we consider the operator 
\begin{eqnarray*}
K&=&(\alpha I+(N-\Nbar))^{-1}(\alpha I+N)\\
&=&(\alpha I+(N-\Nbar))^{-1}(\alpha I+(N-\Nbar)+\Nbar)\\
&=&(\alpha I+(N-\Nbar)^{-1}(\alpha I+(N-\Nbar))+(\alpha I+(N-\Nbar))^{-1}\Nbar\\
&=&I+(\alpha I+(N-\Nbar))^{-1}\Nbar.
\end{eqnarray*}
$K$ is invertible as the product of two invertible operators in $\NovHt$. 
By the ideal property of $\NivH$, Proposition \ref{5.4.5}, 
$(\alpha I+(N-\Nbar))^{-1}\Nbar\in\NivH$; therefore Theorem \ref{Baskakov} implies that $K^{-1}\in \NivHt$. 
The composition of $K^{-1}\in\NivHt$ and $(\alpha I+(N-\Nbar))^{-1}\in\NovHt$ is also in $\NovHt$. 
Therefore, $(\alpha I+N)^{-1}=K^{-1}(\alpha I+(N-\Nbar))^{-1}\in\NovHt$.  
\end{proof}

Theorem \ref{KurbatovH} allows us to prove the spectral algebra property for convolution operators on the Heisenberg group.

\begin{corollary}\label{thmH}
Let $\H$ be the general, reduced Heisenberg group, $v$ an admissible weight function, 
and $S_{f}$ the operator given by convolution with $f$.
If $\alpha_{1}I+S_{f}$, $f\in\LovH$, is invertible in $\B(L^p(\H))$, 
then $(\alpha_{1}I+S_{f})^{-1}=\alpha_{2}I+S_{g}$, $g\in\LovH$. 
\end{corollary}
$\remark$
Barnes proves in \cite{Bar90} that the spectral algebra property for a convolution operator on 
$L^{1}(G)$ is equivalent to $G$ being amenable and symmetric. 
Since $\H$ is nilpotent, it is symmetric \cite{Lud79}. 
Taking $M$ as a mean on $L^{\infty}(\G\times\Gh)$, $M_{\H}(f)=\int_{\T}M(f(\cdot,\cdot,e^{2\pi i\tau}))d\tau$, 
is a shift-invariant mean on $L^{\infty}(\H)$, and consequently $\H$ is amenable; see chapters 2 and 12 of \cite{Pie84}.
Therefore, Corollary \ref{thmH} also follows from  Barnes's work in \cite{Bar90}.
For the case when $G$ is compactly generated, Corollary \ref{thmH} is also a special case of Theorems 
3.6 and 3.7 in \cite{FGLLM03}. 
\vspace{.05in}

\begin{proof} 
For $F_{1},F_{2}\in \LovH$, $\|F_{1}\hstar F_{2}\|_{\LovH}
\leq\int_{\H}\int_{\H}|F_{1}(\h)||F_{2}(\h^{-1}\hn)|v(\hn)d\hn d\h=
\int_{\H}|F_{1}|(\int_{\H}|F_{2}(\hn)|v(\hn)d\hn)v(\h)d\h
=\|F_{1}\|_{\LovH}\|F_{2}\|_{\LovH}$.
Consequently, 
$$(\ao\delta+F_{1})\hstar(\at\delta+F_{2})=\ao\at\delta+\ao F_{2}+\at
F_{1}+F_{1}\hstar F_{2}=\ath\delta+F, \quad F\in L^{1}_{v}(\H).$$
To meet the conditions of Theorem \ref{KurbatovH}, we define the function 
$F(\hn,\h)=f(\hn\h^{-1})$.  
Then $\alpha I$ plus the integral operator with kernel $F$ is the same as $S_{\alpha\delta +f}$ 
and satisfies the conditions of Theorem \ref{KurbatovH}. 
Assuming $S_{\ao\delta+f}$ to be invertible in $\B(L^{p}(H))$, 
Theorem \ref{KurbatovH} states that $S^{-1}_{\ao\delta+f}=\at I+A$, where $A$ is an integral 
operator majorized by a function $\beta\in \LovH$.
We use an approximate identity $\{\psi_{n}\}_{n\geq0}\subset \LovH$.
We set $\theta_{n}=S^{-1}_{\ao\delta+f}\psi_{n}$, and
$\theta=\lim_{n\rightarrow \infty}(\at I+A)\psi_{n}=\at\delta+\lim_{n\rightarrow\infty}A\psi_{n}$. 
Since $A$ is majorized by $\beta\in \LovH$, 
\begin{eqnarray*}
\lim_{n\rightarrow\infty}|A\psi_{n}|(\hn)&\leq&\lim_{n\rightarrow\infty}\int_{\H}\psi_{n}(\h)\beta(\h^{-1}\hn)d\h\\
&=&\beta(\hn)
\end{eqnarray*}
Set $g=\lim_{n\rightarrow\infty}A\psi_{n}\in L^{1}_{v}(\H)$. 
Then $\theta=\at\delta+g$, and by the continuity of convolution, $(\ao\delta+f)\hstar (\at\delta+g)=\delta$. 
For any $\phi\in C_{0}(\H)$, the space of continuous, compactly supported functions on $\H$, 
\begin{eqnarray*}
S_{\ao\delta+f}(S_{\theta}-(\at I+A))\phi&=&S_{\at\delta+f}S_{\at\delta+g}\phi-S_{\ao\delta+f}S^{-1}_{\ao\delta+f}\phi\\
&=&(\ao\delta+f)\hstar(\at\delta+g)-\phi\\
&=&\delta\hstar\phi-\phi\\
&=&0
\end{eqnarray*}
$S_{\ao\delta+f}$ is assumed invertible in $\B(\LpH)$, and both $S_{\at\delta+g}\phi$ and $(\at I+A)\phi$
are in $\LovH$; therefore,$(S_{\at\delta+g}-(\at I+A))\phi=0$ for all $\phi\in C_{0}(\H)$. 
Since the space of continuous compactly supported functions is dense in $\LovH$, 
$S_{\at\delta+g}=\at I+A$, and $S^{-1}_{\ao\delta+f}=S_{\at\delta+g}$, $g\in \LovH$. 
Equivalently, if $\ao I+S_{f}$, $f\in \LovH$, is invertible in $\B(\LpH)$, 
then its inverse is also of the form $\at I+S_{g}$, and $g\in \LovH$.
\end{proof}
\end{section}

\begin{section}{Spectral Algebra Property for Twisted Convolution and Pseudodifferential Operators}                           
$\T$ is originally adjoined to $\GGh$, thus creating the Heisenberg group, in order to obtain group structure for $\GGh$. 
However, functions defined only $\GGh$ are still of special interest, particularly for pseudodifferential operators. 
Here we discuss the Weyl pseudodifferential operator $\Ls$, given by a symbol $\sigma\in\mathcal{S}'(\GGh)$:
\begin{equation}\label{WeylOp}
\Ls f(t)=\int_{\G}\int_{\Gh}\hat{ \sigma}(\omega,x)e^{-\pi i \omega\cdot x}T_{-x}M_{\omega}f(t)d\omega dx.
\end{equation} 
The map $\sigma \mapsto \Ls$ is called the Weyl transform, and $\sigma$ and $\s$ are called the 
symbol and spreading function of the operator $\Ls$. 
The composition rule for two Weyl pseudodifferential operators 
is $\Ls L_{\tau}=L_{\mathcal{F}^{-1}(\s\natural\t)}$, 
where $\natural$ denotes \emph{twisted convolution}\cite{Fol89} and is defined by
\begin{equation*}
 F\natural G(x_{0},\omega_{0})= \int_{\G}\int_{\Gh}F(x,\omega)G(x_{0}-x,\omega_{0}-\omega)e^{\pi i(x\omega_{0}-\omega x_{0})}d\omega dx.
\end{equation*}
Since $F\natural G(x,\omega)\leq |F|*|G|(x,\omega)$, 
twisted convolution is dominated by regular convolution.
Therefore, $\LovGGh$ is closed with respect to twisted convolution.
In order to prove the spectral algebra property for twisted convolution on $\LovGGh$
we need a weighted version of Kurbatov's Theorem \ref{Kurbatov}.
\begin{theorem}\label{KurbatovG}
Let $\G$ be a locally compact abelian group. Then $\NovGt$ is an inverse-closed Banach algebra in $\B(L^{p}(\G),l^{q})$, 
$1\leq p,q\leq\infty$.
\end{theorem}
\begin{proof}
Kurbatov's proof of his Theorem 5.4.7 \cite{Kur99} holds here. The addition of weights is justified by Theorem \ref{BaskakovFull}.
\end{proof}
\begin{corollary}\label{wienerGG}
Let $\G$ be a locally compact abelian group and $\Gh$ its dual group, 
and let $T_{f}\in\B(L^{p}(\GGh))$ be the operator given by twisted convolution with $f$: 
$T_{f}\phi=f\natural \phi$, $\phi \in L^{p}(\GGh)$. 
If $\alpha_{1}I+T_{f}$ is invertible in $\B(L^{p}(\GGh))$ and $f\in\LovGGh$, 
then $(\alpha_{1}I+T_{f})^{-1}=\alpha_{2}I+T_{g}$ and $g\in\LovGGh$. 
\end{corollary}
\begin{proof}
The proof of Corollary \ref{thmH} carries over exactly with the sole substitution of $\GGh$ and $\natural$ for $\H$ and $\hstar$.
\end{proof}

Before applying these theorems to pseudodifferential operators, we briefly discuss the importance of 
pseudodifferential operators in the study of time-varying communication systems, such as wireless communications. 
We view $f(t)$ as a transmitted signal; then $T_{x}f(t)$, $x>0$, corresponds to a time shift of the signal, 
and $M_{\omega}f(t)$ corresponds to a modulation or frequency shift. 
The received signal at time $t_{0}$ is a weighted collection of delayed, modulated copies of the transmitted signal. 
Therefore the received signal may be expressed as 
\begin{equation}
f_{rec}(t_{0})=\int_{\R}\int_{\R}\s(x,\omega)T_{-x}M_{\omega}f_{trans}(t_{0})dxd\omega\label{channel2},
\end{equation}
where we have absorbed $e^{-\pi i \omega\cdot x}$ into $\s$. 
The assumption that $\s\in L^{1}_{v}(\R^{2})$ is appropriate, as in practice the strength of the delayed 
copies of the signal decays quickly in time. 
The Doppler effect or frequency shift depends on the travel speed of the signal and the relative speeds and angles between the 
transmitter, any reflecting bodies, and the receiver. 
Since these quantities are all bounded in practice, the Doppler effect is also bounded. 
Hence if the Doppler effect is limited to $[-D,D]$, the support of 
$\s(x,\cdot)$ is contained in $[-D,D]$ for all $x$ \cite{Rap96,Str06}. 
In practice one must ``numerically invert'' the operator in (\ref{channel2}). 
Theorem \ref{thm2} states that the inverse will have the same off-diagonal decay as the original operator. 
The resulting matrix may, therefore, be truncated to a small number of diagonals, which is essential for 
fast real-world computation.

Above we showed that $\widetilde{L^{1}_{v}}(\GGh)$ is an inverse-closed Banach algebra with respect to twisted convolution. 
Due to the composition rule, the previous theorems easily establish the spectral algebra property for a class of 
Weyl pseudodifferential operators.

\begin{theorem} \label{thm2}
Let $\textnormal{OP}(\mathcal{F}^{-1}L^{1}_{v}(\GhG))$ 
denote the space of pseudodifferential operators with Weyl symbol $\sigma$ satisfying 
$\s\in L^{1}_{v}(\GhG)$.
Then $\widetilde{\textnormal{OP}}(\mathcal{F}^{-1}L^{1}_{v}(\GhG))$ is an inverse-closed subalgebra of $\B(L^{p}(\G))$. 
That is
\begin{enumerate}
\item[(i)]
$\alpha I+\Ls$ is bounded on all $L^{p}(\G)$.
\item[(ii)]
If $\s ,\t\in L^{1}(\GhG)$, then $(\ao I+\Ls)(\at I+L_{\tau})=(\ath I+L_{\gamma})$, where $\g\in L^{1}(\GhG)$.
\item[(iii)]
If $\ao I+\Ls$ is invertible in $\B(L^{p}(\G))$, then $(\ao I+\Ls)^{-1}=(\at I+L_{\tau})$ where $\t\in L^{1}(\GhG)$.
\end{enumerate}
\end{theorem}

\begin{proof}

\noindent (i).\begin{eqnarray*}\|\Ls f\|^{p}_{L^{p}}&\leq&\int_{\G}\left|\int_{\Gh}\int_{\G}\s(\omega,x)e^{-\pi i \omega\cdot x}T_{-x}M_{\omega}f(t)dxd\omega\right|^{p}dt\\
&\leq& \int_{\G} \left(\int_{\G} \int_{\Gh} |\s(\omega,x)||f(t+x)|d\omega dx \right)^{p}dt\\
&=&\int_{\G} \left( \int_{\G} \|\s(\cdot,-x)\|_{L^{1}}|f(t+x)|dx\right)^{p}dt\\
&=&\|\;\|\s(\cdot,u)\|_{L^{1}}*|f|(u)\|^{p}_{L^{p}}\\
&\leq&\|\s\|^{p}_{L^{1}}\|f\|^{p}_{L^{p}}\\
&\leq&\|\s\|^{p}_{L^{1}_{v}}\|f\|^{p}_{L^{p}}
\end{eqnarray*}
Therefore, $\|(\alpha I+\Ls)f\|_{L^{p}(\G)}\leq (|\alpha|+\|\s\|_{L^{1}_{v}(\G)})\|f\|_{L^{p}(\G)}$.

\noindent (ii). 
By Corollary \ref{wienerGG}, if $\s ,\t\in L^{1}(\GhG)$, 
$(\ao\delta+\s)\natural(\at\delta+\t)=(\ath \delta+\hat{\gamma})$, 
where $\g \in L^{1}(\GhG)$. 
Then $\mathcal{F}^{-1}(\alpha_{3} \delta+\g)=\alpha_{3}+\gamma$, and  
$(\ao I+\Ls)(\at I+L_{\tau})=(\ath I+L_{\gamma})$.

\noindent (iii). 
Let $(\ao I+\Ls)^{-1}=A$, $A\in \B(L^{p}(\G))$. Using the Schwartz kernel theorem, Gr\"ochenig shows in \cite{Gro01}, 
that there exists a symbol $\gamma\in\mathcal{S}'(\G)$, such that $A=L_{\gamma}$.
In order to apply Corollary \ref{wienerGG} we must show that the twisted convolution operator $T_{\hat{\gamma}}$ 
is bounded as an operator on $L^{1}(\Gh\times\G)$. 
By the closed graph theorem \cite{Rud91}, $T_{\hat{\gamma}}$ is \emph{not} bounded on $L^{1}(\Gh\times\G)$ 
if and only if there exists a sequence 
$\{\hat{\phi}_{n}\}_{n\in\Nat}\subset L^{1}(\Gh\times\G)$ such that $\{\hat{\phi}_{n}\}\rightarrow 0$, 
but $T_{\hat{\gamma}}\hat{\phi}_{n}\nrightarrow 0$. 
If $T_{\hat{\gamma}}\notin \B(L^{1}(\Gh\times\G))$, then there exists an $\epsilon >0$ and 
a subsequence $\{\hat{\phi}_{n_{k}}\}$ such that $\|T_{\hat{\gamma}}\hat{\phi}_{n_{k}}\|_{ L^{1}(\Gh\times\G)}>\epsilon$
for all $n_{k}$.
$L_{\gamma}\in\B(L^{p}(\G))$ by assumption, and $L_{\phi_{n_{k}}}\in\B(L^{p}(\G))$ by (i). 
Therefore $L_{\gamma}L_{\phi_{n_{k}}}\in \B(L^{p}(\G))$. 
By (i) we have
\begin{eqnarray*}
\| L_{\gamma}L_{\phi_{n_{k}}}-L_{\gamma}L_{\phi_{n_{l}}}\|_{\B(L^{p}(\G))}\leq\|L_{\gamma}\|_{\B(L^{p}(\G))}\|\hat{\phi}_{n_{k}}-\hat{\phi}_{n_{k}}\|_{ L^{1}(\Gh\times\G)}.
\end{eqnarray*}
Since $\{\hat{\phi}_{n_{k}}\}$ is a convergent sequence, $\{ L_{\gamma}L_{\phi_{n_{k}}}\}$ is a Cauchy sequence 
in the Banach space $ \B(L^{p}(\G))$. 
Since $\lim_{k\rightarrow \infty}\{\hat{\phi}_{n_{k}}\}\rightarrow 0$, 
$\lim_{k\rightarrow \infty}\{ L_{\gamma}L_{\phi_{n_{k}}}\}\equiv 0$, 
where the latter convergence is in operator norm. 
By Theorem 14.6.1 in \cite{Gro01}, $\lim_{k\rightarrow \infty}\|T_{\hat{\gamma}}\hat{\phi}_{n_{k}} \|_{ L^{2}(\Gh\times\G)}=0$, 
which implies  $\lim_{k\rightarrow \infty}\|T_{\hat{\gamma}}\hat{\phi}_{n_{k}} \|_{ L^{1}(\Gh\times\G)}=0$
and contradicts the assumption. Thus $T_{\hat{\gamma}}\in \B(L^{1}(\Gh\times\G))$.
Therefore $(\ao I+\Ls)L_{\gamma}=I$ implies $(\ao\delta+\s)\natural\hat{\gamma}=\delta$, and similarly 
$\hat{\gamma}\natural(\ao\delta+\s)=\delta$. 
By Corollary 3.2., $\gamma=\at\delta+\tau$, $\t\in L^{1}(\GhG)$.

\end{proof}
\end{section}

\noindent \emph{Acknowledgements}: B.~F.\ was
supported by NSF VIGRE grant DMS-0135345. 
T.~S.\ was supported by NSF grant DMS-0511461.


\noindent DEPARTMENT OF MATHEMATICS\\
UNIVERSITY OF CALIFORNIA\\
DAVIS, CA 95616, USA. 

\noindent \emph{email}: farrell@math.ucdavis.edu, strohmer@math.ucdavis.edu

\end{document}